\documentclass[12pt,a4paper,reqno]{amsart}

\usepackage{amssymb,mathtools}

\usepackage[marginratio=1:1]{geometry}

\edef\restoreparindent{\parindent=\the\parindent\relax}
\usepackage{parskip}
\restoreparindent

\makeatletter
\@namedef{subjclassname@2020}{\textup{2020} Mathematics Subject Classification}
\makeatother

\newtheorem{thm}{Theorem}
\newtheorem{lem}{Lemma}

\newtheorem{cor}{Corollary}
\newtheorem{conj}{Conjecture}

\newtheorem{Thm}{Theorem}

\newtheorem{Lem}{Lemma}

\theoremstyle{definition}
\newtheorem{rem}{Remark}
\newtheorem{defn}{Definition}

\newcommand{\IC}{{\mathbb C}}
\newcommand{\ID}{{\mathbb D}}

\newcommand{\IT}{{\mathbb T}}
\newcommand{\IS}{{\mathcal S}}

\newcommand{\real}{{\operatorname{Re}\,}}

\newcommand{\eit}{{e^{i\theta}}}

\newcommand{\be}{\begin{equation}}
\newcommand{\ee}{\end{equation}}

\newcommand{\blem}{\begin{lem}}
\newcommand{\elem}{\end{lem}}

\newcommand{\bdefn}{\begin{defn}}
\newcommand{\edefn}{\end{defn}}

\newcommand{\bthm}{\begin{thm}}
\newcommand{\ethm}{\end{thm}}

\newcommand{\bcor}{\begin{cor}}
\newcommand{\ecor}{\end{cor}}

\newcommand{\bconj}{\begin{conj}}
\newcommand{\econj}{\end{conj}}

\newcommand{\brem}{\begin{rem}}
\newcommand{\erem}{\end{rem}}

\newcommand{\bpf}{\begin{proof}}
\newcommand{\epf}{\end{proof}}

\begin{document}
	
\bibliographystyle{abbrv}

\title[On harmonic quasiregular mappings in Bergman spaces]
{On harmonic quasiregular mappings in Bergman spaces}

\author{Suman Das}
\address{Suman Das\vskip0.05cm Department of Mathematics with Computer Science, Guangdong Technion - Israel
	Institute of Technology, Shantou, Guangdong 515063, P. R. China.}
\email{suman.das@gtiit.edu.cn}

\author{Antti Rasila}
\address{Antti Rasila \vskip0.05cm Department of Mathematics with Computer Science, Guangdong Technion - Israel
	Institute of Technology, Shantou, Guangdong 515063, P. R. China, \vskip0.01cm and \vskip0.01cm Department of Mathematics, Technion - Israel
	Institute of Technology, Haifa 3200003, Israel.}
\email{antti.rasila@gtiit.edu.cn; antti.rasila@iki.fi}



\begin{abstract}
A classical result of Hardy and Littlewood says that if $f=u+iv$ is analytic in the unit disk $\ID$ and $u$ is in the harmonic Bergman space $a^p$ ($0<p<\infty$), then $v$ is also in $a^p$. This complements a celebrated result of M.~Riesz on Hardy spaces, which only holds for $1<p<\infty$. These results do not extend directly to complex-valued harmonic functions. We prove that the Hardy-Littlewood theorem holds for a harmonic function $f=u+iv$ if we place the assumption that $f$ is quasiregular in $\ID$. This makes further progress on the recent Riesz type theorems for harmonic quasiregular mappings by several authors.

Then we consider univalent harmonic mappings in $\ID$ and study their membership in Bergman spaces. In particular, we produce a non-trivial range of $p>0$ such that every univalent harmonic function $f$ (and the partial derivatives $f_\theta,\, rf_r$) is of class $a^p$. This result extends nicely to harmonic quasiconformal mappings in $\ID$.
\end{abstract}

\maketitle
\pagestyle{myheadings}
\markboth{S. Das and A. Rasila}{On harmonic quasiregular mappings in Bergman spaces}

%%%%%%%%%%%%%%%%%%%%%%%%%%%%%%%%%%%%%%%%%%%%%%%%%%%%%%%%%%%%%%%%%%%%%%%%%%%%%%%%%%%%%%%%

%\tableofcontents

\section{Introduction and Background}\label{sec1}

\subsection{Hardy and Bergman Spaces}
Suppose $\ID$ denotes the open unit disk in the complex plane and $\IT$ is the unit circle. If $f$ is analytic in $\ID$,
the \textit{integral means} of $f$ are defined as $$M_p(r, f) \coloneqq \left( \frac{1}{2\pi}\int_{0}^{2\pi}\vert f(r e^{i\theta}) \vert ^p\, d\theta \right)^{1/p}, \quad 0 < p < \infty.$$ %and $$M_\infty(r,f) \coloneqq \sup_{ \vert z \vert=r} \vert f(z) \vert.$$
The function $f$ is in the \textit{Hardy space} $H^p$ if $$\|f\|_{H^p} \coloneqq \lim_{r \to 1^-} M_p(r,f)<\infty.$$ Similarly, $f$ is said to be in the \textit{Bergman space} $A^p$ ($0<p<\infty$) if $$\|f\|_p \coloneqq \left(\int_{\ID}|f(z)|^p\, dA(z)\right)^{1/p}<\infty,$$ where $$dA(z)=\frac{1}{\pi}dx\,dy=\frac{1}{\pi}r\,dr\,d\theta \quad (z=x+iy=r\eit \in \ID)$$ is the normalized area measure. Since $$\|f\|_p^p=2\int_{0}^1 M_p^p(r,f)\, r\, dr,$$ it is obvious that $H^p \subset A^p$. In fact, it is even true that $H^p\subset A^{2p}$ with \be \label{Hardy_Bergman}\|f\|_{2p} \le \|f\|_{H^p}, \quad 0<p<\infty.\ee
One can find a detailed survey on Hardy spaces in the book of Duren \cite{Duren}. We refer to the books of Duren \cite{Duren_Bergman} and Hedenmalm, Korenblum, and Zhu \cite{HKZbook} for the theory of Bergman spaces. In this paper, we mostly follow notations from \cite{Duren_Bergman} and \cite{Duren}.

%A real-valued function $u(x,y)$, twice continuously differentiable in $\ID$, is called \textit{harmonic} if it satisfies the Laplace equation $$\Delta u =\frac{\partial^2 u}{\partial x^2}+\frac{\partial^2 u}{\partial y^2}= 0 $$ in $\ID$. %For a complex-valued function $f=u+iv$, the complex partial derivatives have the form $$f_z=\frac{1}{2} \left(f_x-if_y\right) \quad \text{and} \quad f_{\overline{z}}=\frac{1}{2} \left(f_x+if_y\right),$$ where $z=x+iy\in \IC$. In view of this, it is easy to see that the Laplacian operator $\Delta$ can also be written as $$\Delta=4\, \frac{\partial^2}{\partial \overline{z} \partial z} \cdot$$
A complex-valued function $f=u+iv$ in $\ID$ is harmonic if $u,v$ are real-valued harmonic functions in $\ID$. Every such function has a unique representation $f=h+\overline{g}$, where $h$ and $g$ are analytic in $\mathbb{D}$ with $g(0)=0$. The function $f$ is said to be in the \textit{harmonic Hardy space} $h^p$ if $M_p(r,f)$ is bounded, and in the \textit{harmonic Bergman space} $a^p$ if $\|f\|_p<\infty$. We note that for harmonic functions, no inequality of type \eqref{Hardy_Bergman} holds for $0<p\le 1$ (see \cite{Kalaj_Mestro}), so $h^p \not \subset a^{2p}$. However, the obvious inclusion $h^p \subset a^p$ holds.
%We refer to the 
%book of Duren \cite{Duren:Harmonic} for the theory of planar harmonic mappings, and to the 
%monograph of Pavlovi\'c \cite{Pavbook} for a concise survey on $h^p$ spaces.

\subsection{Growth of conjugate functions}
Let $u$ be a real-valued harmonic function in $\ID$, and let $v$ denote its harmonic conjugate normalized by $v(0)=0$. It is of natural interest that if $u$ has a certain property, whether so does $v$. In the study of boundary behaviour, this is ensured by a celebrated theorem of M.~Riesz.
\begin{Thm}\label{MRiesz}\cite[Theorem 4.1]{Duren}
	If $u \in h^p$ for some $p$, $1<p<\infty$, then its harmonic conjugate $v$ is also of class $h^p$. Furthermore, there is a constant $A_p$, depending only on $p$, such that $$M_p(r,v) \le A_p \, M_p(r,u),$$ for all $u \in h^p$.
\end{Thm}
The theorem fails for $p=1$ and $p=\infty$, as illustrated by counterexamples given in \cite[p. 56]{Duren}. It was then proved by Kolmogorov that the harmonic conjugate of an $h^1$-function, while not necessarily in $h^1$, does belong to $h^p$ for all $p<1$. Indeed, Zygmund showed that the condition $|u| \log^+ |u| \in L^1(\IT)$ is ``minimal" to guarantee that $v\in h^1$. %We refer to the paper of Pichorides \cite{Pichorides} for the optimal constants in the Riesz, Kolmogorov, and Zygmund theorems.

Hardy and Littlewood \cite{HL2} explained that in the case $0<p<1$, Riesz's theorem is false in a more comprehensive sense. One might expect from Kolmogorov's result that if $u\in h^p$, then $v$ should belong to $h^q$ for $0<q<p$. But this is not true, and in fact, $v$ is not necessarily in $h^q$ for any $q>0$. Interestingly, they showed that the symmetry in Riesz's theorem is restored for all $p$ if one considers the ``orders of infinity" of the means. That is, if $M_p(r,u)=O\left((1-r)^{-\beta}\right)$ for $0<p\le \infty$ and $\beta>0$, then so is $M_p(r,v)$.

It follows immediately from the theorem of Riesz that the harmonic Bergman space $a^p$ is also ``self-conjugate" for $1 < p < \infty$, and $\|v\|_p \le C_p\|u\|_p$ for all $u\in a^p$. However, it is remarkable that in the Bergman norm the theorem remains true even for $0 < p \le 1$. Indeed, this is one of the rare occasions where Bergman spaces behave better than Hardy spaces.

\begin{Thm}\label{conj}\cite[Theorem 5]{HL2}
	Let $0<p< \infty$. If $u$ is a harmonic function of class $a^p$, then its harmonic conjugate $v$ is also in $a^p$, and $$\|v\|_p \le C_p\|u\|_p,$$ where $C_p$ is a constant depending only on $p$.
\end{Thm}

Hardy and Littlewood stated this result in \cite{HL2} and indicated a proof for $0<p\le 1$. The technical details were later given by Watanabe \cite{Watanabe}, and subsequently refined by Sobolewski \cite{Sobolewski}. Another proof, in a somewhat different setting, can be found in the book \cite{Pavlovic_New} of Pavlovi\'c (Theorem 3.8).

%They also proved the following converse.
%\begin{Thm}\label{hl_new2}\cite[Theorems 3 \& 6]{HL2}
%	Let $f=u+iv$ be analytic in $\ID$, and $$M_p(r,f') = O \left(\frac{1}{(1-r)^{\beta+1}}\right), \quad 0<p\le \infty,\quad \beta> 0.$$ Then $$M_p(r,f) = O \left(\frac{1}{(1-r)^{\beta}}\right).$$ Furthermore, if $$M_p(r,f')=O \left(\frac{1}{1-r}\right),$$ then $$M_p(r,f)=O \left(\left(\log\frac{1}{1-r}\right)^\gamma\right),$$ where $\gamma=\max\{1/p,1\}$.
%\end{Thm}
%Let us note that the functions $|u|^p$ and $|v|^p$ are subharmonic when $p\ge 1$, but not when $p<1$, and therefore, $M_p(r,u)$ and $M_p(r,v)$ are not necessarily monotonic for $p<1$. This is the principal difficulty in dealing with the case $0<p<1$ for harmonic functions.

%\subsection{Smoothness class}

\subsection{Recent developments on harmonic quasiregular mappings}
For $K\ge 1$, a sense-preserving harmonic function $f$ is said to be $K$-\textit{quasiregular} if the \textit{dilatation} $$D_f\coloneqq \frac{|f_z|+|f_{\bar{z}}|}{|f_z|-|f_{\bar{z}}|} \le K$$ throughout $\ID$. Writing $f=h+\overline{g}$, this condition is equivalent to saying that its \textit{analytic dilatation} $\omega \coloneqq g'/h'$ satisfies the inequality $$|\omega(z)| \le k <1 \quad (z\in \ID),$$ where \be\label{eq1}k\coloneqq \frac{K-1}{K+1}.\ee The function $f$ is called $K$-\textit{quasiconformal} if it is $K$-quasiregular and homeomorphic in $\ID$. One can find the $H^p$-theory for quasiconformal mappings in, for example, the paper of Astala and Koskela \cite{ast_kos}. %It is worth mentioning that harmonic quasiconformal mappings have generated considerable interest in recent times, perhaps from a novel point of view. In \cite{wang_rasila}, Wang et al. constructed independent extremal functions for harmonic quasiconformal mappings, which were then further explored by Li and Ponnusamy in \cite{Li_Pon}. Later in \cite{DHR2}, it was established that the newly constructed functions are also extremal for Baernstein type results in the Hardy space of harmonic quasiconformal mappings.

If $f=u+iv$ is a harmonic function in $\ID$, and $u\in h^p$ for some $p>1$, the imaginary part $v$ does not necessarily belong to $h^p$. In other words, the classical Riesz theorem does not extend directly to harmonic functions. This raises a natural question: \textit{Under what additional assumptions does a harmonic version of the Riesz theorem hold?} In the insightful paper \cite{Liu_Zhu}, Liu and Zhu identified the quasiregularity of $f$ as a sufficient condition for this. The result, sharpened significantly by Chen and Huang \cite{chen_huang}, can be stated as follows.
\begin{Thm}\label{LZ}\cite{Liu_Zhu}
	Let $f=u+iv$ be a harmonic $K$-quasiregular mapping in $\ID$ with $v(0)=0$. If $u \in h^p$ for some $p\in (1,\infty)$, then $v$ is also in $h^p$, and the inequality $$M_p(r,v) \le C_{p,K} M_p(r,u)$$ holds for some constant $C_{p,K}$ that depends only on $p$ and $K$.
\end{Thm}
%The assumption that $f$ is harmonic cannot be dropped. Indeed, there exists a quasiconformal mapping for which the first coordinate function is in the Hardy space, but the function itself is not (see \cite[p. 36]{ast_kos}).
Later in \cite{Kalaj_Kolmogorov}, Kalaj produced a couple of Kolmogorov type theorems for harmonic quasiregular mappings. Subsequently, a quasiregular analogue of Zygmund's theorem was obtained by Kalaj \cite{Kal_Zyg}, and also independently by the present authors and Huang \cite{DHR1}. In the recent paper \cite{DR1}, the present authors  proved that the real and imaginary parts of a harmonic quasiregular mapping have the same order of infinity for all $p\in (0,\infty]$, thereby extending Theorem \ref{LZ} to the cases $0<p\le1$ and $p=\infty$.
%In \cite{DHR2}, Baernstein type extremal results have been obtained for the integral means of harmonic quasiregular mappings.

We present the main results of this paper and related discussions in the next section. The proofs and supporting results are given in Section \ref{sec_proofs}.

\section{Main Results}\label{results}
The purpose of this paper is twofold. First, we establish the quasiregular analogue of Theorem \ref{conj}, i.e., we show that the real and imaginary parts of a harmonic quasiregular mapping necessarily belong to the same Bergman space $a^p$ for any $p>0$. In view of Theorem \ref{LZ}, the case $0<p\le 1$ is particularly interesting and is the primary motivation of this paper.
\bthm\label{main}
Suppose $0<p<\infty$ and $f=u+iv$ is a harmonic $K$-quasiregular mapping in $\ID$. If $u$ is in the harmonic Bergman space $a^p$, then $v$ is also in $a^p$, and $$\|v\|_p \le C_{p,K} \|u\|_p,$$ where $C_{p,K}$ is a constant that depends only on $p$ and $K$.
\ethm

As an application of Theorem \ref{main}, we obtain the following integral mean estimate for harmonic quasiregular mappings.
\bcor\label{cor1}
Let $f=u+iv$ be a harmonic $K$-quasiregular mapping in $\ID$ and $u\in a^p$ for some $p\in (0,1]$. Then $$M_p(r,f) \le \frac{C_{p,K}}{(1-r)^{1/p}}\|u\|_p,$$ for some constant $C_{p,K}$ depending on $p$ and $K$.
\ecor

In the second part of this paper, we study the membership of univalent harmonic functions and harmonic quasiconformal mappings in Bergman spaces. Suppose $\IS_H$ is the class of all sense-preserving univalent harmonic functions $f=h+\overline{g}$ in $\ID$ with the normalizations  $h(0)=g(0)=h'(0)-1=0$. The quantity $$\alpha \coloneqq \sup_{f\in \IS_H} \left|\frac{h''(0)}{2}\right|$$ is known as the \textit{order} of the class $\IS_H$. In \cite{Abu_Lyzzaik}, Abu-Muhanna and Lyzzaik initiated the study of boundary behaviour of functions in $\IS_H$, and showed that every $f \in \IS_H$ is of class $h^p$ for $p<1/(2\alpha+2)^2$. This range of $p$ was notably improved to $p<1/\alpha^2$ by Nowak \cite{Nowak}, who conjectured the sharp range $p<1/\alpha$. The conjecture was recently verified (under an additional condition) by the first author and Sairam Kaliraj \cite{DK3}.

As the other main result, we show that the class $\IS_H$ is contained in the harmonic Bergman space $a^p$ for $p<1/\alpha$. Indeed, this range of $p$ is new and cannot be deduced from the earlier progress on Hardy spaces (see Remark \ref{rem1}).

\bthm\label{main3}
If $f\in \IS_H$, then the following statements hold:
\begin{enumerate}
	\item[(i)] $f$ is in the harmonic Bergman space $a^p$ for all $p<1/\alpha$.
	\item[(ii)] The analytic functions $f_z, \overline{f_{\bar{z}}} \in A^q$ and the harmonic functions $f_\theta, r f_r \in a^q$ for all $q<1/(\alpha+1)$.	
\end{enumerate}	
\ethm

\brem\label{rem1}
From the result of Nowak, one sees that $\IS_H \subset h^p \subset a^p$ for $p<1/\alpha^2$. It is known that $\alpha\ge 3$, so clearly the range of $p$ in Theorem \ref{main3} is larger. Curiously, our result ensures that the partial derivatives are also of class $a^q$ for some $q>0$. This outcome is typically not possible in the Hardy space setting, since the derivative of a function $f \in H^p$ need not belong to any $H^q$ space. This is another instance where Bergman spaces produce more complete results compared to Hardy spaces.
\erem

For $K\ge 1$, let $\IS_H(K)$ be the class of harmonic quasiconformal mappings in $\ID$, i.e., $$\IS_H(K)\coloneqq \left\{f=h+\bar{g}\in \IS_H: f \text{ is } K \text{-quasiconformal}\right\}.$$ The order of $\IS_H(K)$ is defined as $$\alpha_K \coloneqq \sup_{f\in \IS_H(K)} \left|\frac{h''(0)}{2}\right|.$$ Let us mention that harmonic quasiconformal mappings have received a somewhat novel treatment quite recently. For example, in \cite{wang_rasila}, Wang et al. constructed independent extremal functions for the class $\IS_H(K)$, which were then explicitly studied by Li and Ponnusamy in \cite{Li_Pon}.

%Later in \cite{DHR2}, it was established that the newly constructed functions are also extremal for Baernstein type results in the Hardy space of harmonic quasiconformal mappings.
Similar to Theorem \ref{main3}, we have the following result for the class $\IS_H(K)$.

\bthm\label{main2}
Suppose $f\in S_H(K)$ for some $K \ge 1$. Then
\begin{enumerate}
	\item[(i)] $f\in a^p$ for all $p<1/\alpha_K$,
	\item[(ii)] $f_z, \overline{f_{\bar{z}}} \in A^q$ and $f_\theta, r f_r \in a^q$ for all $q<1/(\alpha_K+1)$.	
\end{enumerate}	
\ethm

\brem\label{rem2}
It is interesting to compare Theorem \ref{main2} with a very well-known result of Astala and Koskela. Theorem 3.2 of \cite{ast_kos} implies that every $f\in \IS_H(K)$ is in $h^p\subset a^p$ for $p<1/(2K)$. However, $1/(2K)$ approaches $0$ for large $K$, whereas $$\frac{1}{\alpha_K} \ge \frac{1}{\alpha}\quad (\text{since }\alpha_K \le \alpha).$$ That is, as $K$ increases, Theorem \ref{main2} eventually gives the better estimate.
\erem

\section{Proofs and Auxiliary Results}\label{sec_proofs}

In what follows, we denote the constants by $C$, $C_p$, $C_{p,K}$, etc., and they may be different from one occurrence to another. Also, throughout the paper we assume that $K$ and $k$ are related by the identity \eqref{eq1}.

Our proof of Theorem \ref{main} is fundamentally based on the following result of Fefferman and Stein \cite{Fef_Ste}, which is, in fact, a formal consequence of Theorem \ref{conj}.

\begin{Thm}\label{Fef_Ste}\cite{Fef_Ste}
Let $u$ be a harmonic function in a domain $\Omega\subset \IC$ and suppose the disk $D_r(a)\coloneqq\{z\in \IC: |z-a|<r\}$ is contained in $\Omega$. Then, for every $p>0$, \be\label{eq_Fef}
|u(a)|^p \le \frac{C_p}{|D_r(a)|}\int_{D_r(a)}|u|^p\, dA,
\ee
where $|D_r(a)|$ denotes the area of $D_r(a)$ and $C_p$ is a constant depending only on $p$.
\end{Thm}

Clearly, if $p\ge 1$, \eqref{eq_Fef} holds with $C_p=1$. Indeed, $|u|$ can be replaced by any non-negative subharmonic function, see \cite[Theorem 3.1]{Pavlovic_New}. For our purpose, we shall also require the following gradient estimate for harmonic functions.
\begin{Thm}\label{grad}\cite[p. 17]{Pavlovic_New}
Suppose $u$ is a real-valued harmonic function in $\ID$, and $a\in \ID$ and $r>0$ are such that $D_r(a) \subset \ID$ . Then there is a constant $K$, independent of $r$ and $a$, such that $$|\nabla u(a)| \le \frac{K}{r} \sup \left\{|u(z)|: z\in D_r(a)\right\}.$$
\end{Thm}

An $n$-dimensional version of this estimate can be found in the book of Axler, Bourdon, and Ramey \cite{Axler_HFT} (p. 42, Exercise 6). Theorems \ref{Fef_Ste} and \ref{grad} lead us to prove the following result for harmonic quasiregular mappings, which will be essential in the proof of Theorem \ref{main}.

\blem\label{lem1}
Let $0<p<\infty$ and $f=u+iv$ be a harmonic $K$-quasiregular mapping in $\ID$. Then, we have \be \label{eq_lem1} |f_z(0)|^p \le C_{p,K} \int_{\varepsilon\ID} |u|^p\, dA \quad (\varepsilon=1/2).
\ee
\elem

\bpf
Let us write $f=h+\overline{g}$ and $F=h+g$, so that $$\real F = \real f=u.$$ We note that $$|F'|=(u_x^2+u_y^2)^{1/2}=|\nabla u|,$$ and therefore, $$|F'(0)|=|\nabla u(0)|\le C \sup_{ \vert z \vert<\varepsilon/2}|u(z)|,$$ by Theorem \ref{grad}, where $\varepsilon=1/2$. Given any $z$ in the disk $|z|<\varepsilon/2$, we apply inequality \eqref{eq_Fef} for the disk $D_{\varepsilon/2}(z) \subset \varepsilon \ID$ to see that $$|u(z)|^p\le C_p \int_{D_{\varepsilon/2}(z)} |u|^p\, dA \le C_p \int_{\varepsilon \ID} |u|^p\, dA.$$ This implies $$|F'(0)|^p \le C_p  \sup_{ \vert z \vert<\varepsilon/2}|u(z)|^p \le C_p \int_{\varepsilon \ID} |u|^p\, dA.$$ Now, we observe that $$F'(z)=h'(z)+g'(z)=h'(z)(1+\omega(z)),$$ so that $$|F'(z)|=|h'(z)||1+\omega(z)|\ge (1-k)|h'(z)|,$$ since $|\omega(z)|\le k$. It follows that $$|h'(0)|^p \le \left(\frac{1}{1-k}\right)^p |F'(0)|^p \le C_{p,K} \int_{\varepsilon\ID} |u|^p\, dA$$ and the proof is complete.
\epf

We are now ready to prove Theorem \ref{main}.

\subsection{Proof of Theorem \ref{main}} For $p>1$, the result immediately follows from Theorem \ref{LZ}, so we focus on the case $0<p\le 1$.

Without loss of generality, we assume that $f(0)=0$. Writing $f=h+\overline{g}$, we start with the inequality \eqref{eq_lem1}, i.e., $$ |h'(0)|^p \le C_{p,K} \int_{\varepsilon\ID} |u(z)|^p\, dA(z) \quad (\varepsilon=1/2).$$ Since the quantity $(1-|z|^2)^{-2}$ is bounded in $\varepsilon\ID$, the above inequality can be rewritten as \be \label{invar} |h'(0)|^p \le C_{p,K} \int_{\varepsilon\ID} |u(z)|^p\, d\tau(z) \quad (\varepsilon=1/2),\ee where $$d\tau(z)=(1-|z|^2)^{-2}\, dA(z).$$ It is known (see \cite[p. 2]{Pavlovic_New}) that the measure $d\tau$ is M\"obius invariant, in the sense that $$\int_{\ID} G\circ \sigma_a \, d\tau = \int_{\ID} G \, d\tau,$$ where $G\ge 0$ is any measurable function in $\ID$ and $$\sigma_a(z)=\frac{a-z}{1-\overline{a}z} \quad (|a|<1)$$ is a conformal automorphism of $\ID$. We mention that $|\sigma_a(z)|$ is known as the \textit{pseudo-hyperbolic distance} between points $a,z\in \ID$. Let us write $$f_{\sigma_a}\coloneqq f\circ \sigma_a = h_{\sigma_a}+\overline{g_{\sigma_a}},$$ where $$h_{\sigma_a}=h\circ \sigma_a, \quad g_{\sigma_a}=g\circ \sigma_a.$$ Then $f_{\sigma_a}$ is harmonic $K$-quasiregular in $\ID$. Applying \eqref{invar} to the function $h_{\sigma_a}$, we find that \be \label{compo} |h'(a)|^p(1-|a|^2)^p \le C_{p,K} \int_{\varepsilon\ID} |u\circ {\sigma_a}(z)|^p\, d\tau(z) = C_{p,K} \int_{B_\varepsilon(a)} |u(z)|^p\, d\tau(z),
\ee
by the invariance property of $d\tau(z)$, where $$B_\varepsilon(a)\coloneqq \sigma_a(\varepsilon \ID) = \{z\in \ID: |\sigma_a(z)|<\varepsilon\},$$ as $\sigma_a(\sigma_a(z))\equiv z$. We Integrate \eqref{compo} over all $a\in \ID$, and change the order of integration, to get \be\label{integ}
\int_{\ID}|h'(a)|^p(1-|a|^2)^p \, dA(a) \le C_{p,K} \int_{\ID} |u(z)|^p\, d\tau(z) \int_{B_\varepsilon(z)} \, dA(a).
\ee
Now, we observe that $B_\varepsilon(z)$ is a \textit{pseudo-hyperbolic disk} centered at $z$ with radius $\varepsilon$. One can see from \cite[Proposition 4.4]{Oper_Zhu} that $B_\varepsilon(z)$ is indeed a Euclidean disk $D_r(z_0)$, where $$z_0=\frac{1-\varepsilon^2}{1-\varepsilon^2|z|^2}z, \quad r=\frac{1-|z|^2}{1-\varepsilon^2|z|^2}\varepsilon.$$ Consequently, its area can be estimated as $$|D_r(z_0)| = \int_{B_\varepsilon(z)} \, dA(a) \le C (1-|z|^2)^2.$$ Thus, it follows from \eqref{integ} that \be\label{integ2}
\int_{\ID}|h'(a)|^p(1-|a|^2)^p \, dA(a) \le C_{p,K} \int_{\ID} |u(z)|^p\, dA(z).
\ee
Now, in order to complete the proof, it is enough to establish the inequality \be \label{final}
\int_{\ID} |f|^p \,dA \le C_{p,K} \int_{\ID}|h'(a)|^p(1-|a|^2)^p \, dA(a).
\ee
This is contained in the following lemma, hence the theorem is proved. \qed

\blem\label{lem2}
Suppose $0<p\le 1$ and $f=h+\overline{g}$ is a harmonic $K$-quasiregular mapping in $\ID$. Then  $$\int_{\ID} |f|^p \,dA \le C_{p,K} \int_{\ID}|h'(a)|^p(1-|a|^2)^p \, dA(a).$$
\elem

\bpf
This is just a slight modification of a well-known result. For $0<p\le 1$, we see that $$\int_{\ID} |f|^p \,dA \le \int_{\ID} |h|^p \,dA +\int_{\ID} |g|^p \,dA.$$ One can find in, for example, \cite[p. 85]{Oper_Zhu} the inequality $$\int_{\ID} |h|^p \,dA \le C_p \int_{\ID}|h'(a)|^p(1-|a|^2)^p \, dA(a)$$ for analytic functions $h$. Similarly, \begin{align*}
\int_{\ID} |g|^p \,dA & \le C_p \int_{\ID}|g'(a)|^p(1-|a|^2)^p \, dA(a)\\ & \le k^p C_p \int_{\ID}|h'(a)|^p(1-|a|^2)^p \, dA(a),
\end{align*} as $|g'(a)|\le k|h'(a)|$. Combining these, we obtain $$\int_{\ID} |f|^p \,dA \le (1+k^p)C_{p} \int_{\ID}|h'(a)|^p(1-|a|^2)^p \, dA(a),$$ as desired.
\epf
 We remark that Lemma \ref{lem2}, although elementary, will be also useful in the proof of Theorem \ref{main2}.
 
\subsection{Proof of Corollary \ref{cor1}} 
We give a proof based on the following inequality of Sobolewski, which is another implication of Theorem \ref{conj}.
\begin{Lem}\label{sobo}\cite{Sobolewski}
If $u\in a^p$, $u(0)=0$, $0<p\le 1$, then $$M_p(r,u) \le \frac{C_p}{(1-r)^{1/p}}\|u\|_p.$$
\end{Lem}

Without loss of generality, we may assume that $f(0)=0$. For $0<p\le 1$, we find that \begin{align*}
M_p^p(r,f)& \le M_p^p(r,u)+M_p^p(r,v)\\ & \le \frac{C_p}{1-r} \left(\|u\|_p^p+\|v\|_p^p\right) \quad (\text{by Lemma } \ref{sobo})\\ & \le  \frac{C_{p,K}}{1-r}\, \|u\|_p^p,
\end{align*} where the last inequality follows from Theorem \ref{main}. The corollary is proved. \qed

The proofs of Theorems \ref{main3} and \ref{main2} are similar, with suitable modifications wherever necessary. We include only the proof of Theorem \ref{main2}, for being somewhat less obvious.

\subsection{Proof of Theorem \ref{main2}}
\textbf{Part (i):} If $f=h+\bar{g} \in S_H(K)$ and any $\zeta \in \ID$, it is easy to check that that the function $$F(z)=\frac{f\left(\frac{z+\zeta}{1+\bar{\zeta}z}\right)-f(\zeta)}{(1-|\zeta|^2)h'(\zeta)} \quad (z \in \ID)$$ also lies in $S_H(K)$. 
%Let us write $$H(z)=z+a_2(\zeta)z^2+a_3(\zeta)z^3+\cdots.$$ 
A simple computation shows that $$F_{z z}(0)=(1-|\zeta|^2)\frac{h''(\zeta)}{h'(\zeta)}-2\bar{\zeta}.$$ As $|F_{z z}(0)| \le 2\alpha_K$, we find that $$\left|\frac{h''(\zeta)}{h'(\zeta)}\right| \le \frac{C}{1-|\zeta|} \quad (\zeta \in \ID),$$ for some constant $C>0$. Since $\zeta \in \ID$ is arbitrary, this can be rewritten as \be\label{condn}
\left|\frac{h''(r\eit)}{h'(r\eit)}\right| \le \frac{C}{1-r}, \quad 0\le r<1. 
\ee
On the other hand, an argument similar to that in \cite[p. 98]{Duren:Harmonic} leads us to the estimate
\be\label{condn2}
\left|h'(r\eit)\right| \le \frac{(1+r)^{\alpha_K-1}}{(1-r)^{\alpha_K+1}}. 
\ee
Now, the Hardy-Stein identity (see \cite[p. 126]{Pomm}) for the function $h'$ implies that $$\frac{d}{dr}\left[r\frac{d}{dr}M_p^p(r,h')\right]=\frac{p^2r}{2\pi} \int_{0}^{2\pi} |h'(r\eit)|^{p-2}|h''(r\eit)|^2\, d\theta \quad (p>0).$$ Since $M_p^p(r,h')$ is a (strictly) increasing function of $r$, we have $$\frac{d}{dr}M_p^p(r,h')>0.$$ Therefore, \begin{align*}
	\frac{d^2}{dr^2}M_p^p(r,h') & \le \frac{p^2}{2\pi} \int_{0}^{2\pi} |h'(r\eit)|^{p-2}|h''(r\eit)|^2\, d\theta\\ & = \frac{p^2}{2\pi} \int_{0}^{2\pi} |h'(r\eit)|^{p}\left|\frac{h''(r\eit)}{h'(r\eit)}\right|^2\, d\theta\\ & \le \frac{p^2}{2\pi} \int_{0}^{2\pi} \frac{(1+r)^{(\alpha_K-1)p}}{(1-r)^{(\alpha_K+1)p}} \cdot \frac{C}{(1-r)^{2}} \, d\theta \quad \left(\text{by } \eqref{condn},\, \eqref{condn2} \right)\\ & \le \frac{C}{(1-r)^{(\alpha_K+1)p+2}}.
\end{align*}
Integrating twice from $0$ to $r$ $(r<1)$, we find that $$M_p^p(r,h') \le \frac{C}{(1-r)^{(\alpha_K+1)p}}.$$ For $0<p\le 1$, an appeal to Lemma \ref{lem2} gives \begin{align*}
\int_{\ID} |f|^p \,dA & \le C_{p,K} \int_{\ID}|h'(z)|^p(1-|z|^2)^p \, dA(z)\\ & \le C_{p,K} \int_{0}^{1} (1-r)^{p} M_p^p(r,h')\, r\, dr\\ & \le C_{p,K} \int_{0}^{1} \frac{dr}{(1-r)^{p \alpha_K}}.
\end{align*} The last integral converges for $p<1/\alpha_K$, and therefore, $f \in a^p$ for $p<1/\alpha_K$.\qed

\textbf{Part (ii):} We start by showing that if $f=h+\bar{g} \in a^p$, then $h, g \in A^p$. This is a simple but useful consequence of Theorem \ref{conj}. 

Let us write $f=u+iv$, and suppose $\tilde{u}$ is the harmonic conjugate of $u$ such that $\tilde{u}(0)=v(0)$, while $\tilde{v}$ is the harmonic conjugate of $v$ with $\tilde{v}(0)=-u(0)$. Let $H=u+i\tilde{u}$ and $G=v+i\tilde{v}$. Then $H$, $G$ are analytic in $\ID$, and \begin{align*}
f=\real H + i\,\real G & = \frac{1}{2} \left(H+\overline{H}\right)+\frac{i}{2} \left( G+\overline{G}\right)\\ & =\frac{1}{2}\left(H+iG\right)+\frac{1}{2}\overline{\left(H-iG\right)}.\end{align*} The above normalizations suggest that $$\frac{1}{2}(H-iG)(0)=0,$$ thus the uniqueness of $h$ and $g$ implies $$h=\frac{1}{2}(H+iG)\quad \text{and} \quad g=\frac{1}{2}(H-iG).$$ As $\vert u\vert, \vert v\vert \le \vert f\vert$, clearly $u,v\in a^p$. In view of Theorem \ref{conj}, we have $\tilde{u},\tilde{v} \in a^p$. It follows that $H$, $G$ are in $A^p$, and therefore, so are $h$ and $g$.

Now, let $f=h+\bar{g}\in S_H(K)$ is in $a^p$. As $h,g\in A^p$, it follows (see \cite[p. 78]{Duren_Bergman}) that $$h',g' \in A^q \quad \text{for } q<\frac{p}{p+1}.$$ We observe that $$q<\frac{1}{\alpha_K+1} \quad \text{since} \quad p<\frac{1}{\alpha_K}.$$ Finally, to conclude that $f_\theta, r f_r \in a^q$, we consider the identities \begin{align*}
-if_\theta & = zh'-\overline{zg'},\\ r f_r & = zh'+\overline{zg'}.\end{align*} It follows that $$\|f_\theta\|_q^q \le \|h'\|_q^q+\|g'\|_q^q <\infty,$$ as $q<1$. Similarly, $\|rf_r\|_q^q < \infty$. The proof of the theorem is now complete.\qed

\subsection*{Acknowledgements} The research was partially supported by the Li Ka Shing Foundation STU-GTIIT Joint Research Grant (Grant no. 2024LKSFG06) and the NSF of Guangdong Province (Grant no. 2024A1515010467).

%%%%%%%%%%%%%%%%%%%%%%%%%%%%%%%%%%%%%%%%%%%%%%%%%%%%%%%%%%%%%%%%%%%%%%%%%%%%%%%%%%%%%%%%

\bibliography{references}

%%%%%%%%%%%%%%%%%%%%%%%%%%%%%%%%%%%%%%%%%%%%%%%%%%%%%%%%%%%%%%%%%%%%%%%%%%%%%%%%%%%%%%%%

\end{document}